\documentclass[11pt]{amsart}
\usepackage{amsmath,amsthm, amscd, amssymb, amsfonts}
\usepackage[all]{xy}
\usepackage{hhline}

\newcommand{\Ind}{\operatorname{Ind}}
\newcommand{\Res}{\operatorname{Res}}

\newcommand\toba{{\mathfrak B }}
\newcommand\trasp{\pi}

\newcommand{\trid}{\triangleright}

\newcommand{\ku}{\mathbb C}

\newcommand{\Z}{{\mathbb Z}}
\newcommand{\N}{{\mathbb N}}

\newcommand{\G}{{\mathbb G}}

\newcommand{\Oc}{{\mathcal O}}

\newcommand\sgn{\operatorname{sgn}}

\theoremstyle{plain}

\newtheorem{lema}{Lemma}[section]
\newtheorem{theorem}[lema]{Theorem}

\newtheorem{prop}[lema]{Proposition}

\theoremstyle{definition}

\theoremstyle{remark}
\newtheorem{obs}[lema]{Remark}
\newtheorem{rmk}[lema]{Remarks}

\newcommand\id{\operatorname{id}}

\newcommand\Id{\operatorname{Id}}

\newcommand\sn{\mathbb S_n}

\newcommand\dn{\mathbb D_n}
\newcommand\an{\mathbb A_n}
\newcommand\am{\mathbb A_m}
\newcommand\at{\mathbb A_3}
\newcommand\ac{\mathbb A_4}
\newcommand\aco{\mathbb A_5}

\newcommand\A{\mathbb A}

\newcommand\sm{\mathbb S_m}

\def\pf{\begin{proof}}
\def\epf{\end{proof}}

\theoremstyle{remark}

\begin{document}

\renewcommand{\baselinestretch}{1.2}

\thispagestyle{empty}

\title[pointed Hopf algebras]
{On pointed Hopf algebras associated with alternating and dihedral
groups}
\author[Andruskiewitsch and Fantino]{ Nicol\'as Andruskiewitsch and Fernando Fantino}
\thanks{This work was partially supported by
Agencia C\'ordoba Ciencia, ANPCyT-Foncyt, CONICET and Secyt (UNC)}
\address{\noindent Facultad de Matem\'atica, Astronom\'\i a y F\'\i sica,
Universidad Nacional de C\'ordoba. CIEM -- CONICET.
\newline \noindent Medina Allende s/n
(5000) Ciudad Universitaria, C\'ordoba, Argentina}
\email{andrus@mate.uncor.edu} \email{fantino@mate.uncor.edu}

\date{\today}

\begin{abstract}
We classify finite-dimensional complex pointed Hopf algebra with
group of group-like elements isomorphic to $\aco$. We show that
any pointed Hopf algebra with infinitesimal braiding associated
with the conjugacy class of $\pi \in \an$ is infinite-dimensional
if the order of $\pi$ is odd except for $\pi=(1 \, 2 \, 3)$ in
$\ac$. We also study pointed Hopf algebras over the dihedral
groups.
\end{abstract}
\maketitle

\centerline{\emph{Dedicado a Mar\'\i a In\'es Platzeck en sus \#
$\aco$ a\~nos}.}

\section*{Introduction}\label {0}

In this article, we continue the work of \cite{AZ,AF} on the
classification of finite-dimensional complex pointed Hopf algebras
$H$ with $G=G(H)$ non-abelian. We follow the Lifting Method -- see
\cite{AS-cambr} for a general reference; in particular, we focus
on the problem of determining when the dimension of the Nichols
algebra associated with conjugacy classes of $G$ is infinite. The
paper is organized as follows. In Section \ref{conventions}, we
review some general facts on Nichols algebras corresponding to
finite groups. We discuss the notion of \emph{absolutely real}
element of a finite group in subsection \ref{absreal}. We then
provide generalizations of \cite[Lemma 1.3]{AZ}, a basic tool in
\cite{AZ,AF}, see Lemmata \ref{lemagen2} and \ref{lemagen1}.
Section 2 is devoted to pointed Hopf algebras with coradical $\ku
\an$. We prove that any finite-dimensional complex pointed Hopf
algebra $H$ with $G(H)\simeq \aco$ is isomorphic to the group
algebra of $\aco$; see Theorem \ref{mainteor}. This is the first
finite non-abelian group $G$ such that all pointed Hopf algebras
$H$ with $G(H)=G$ are known. We also prove that $\dim
\toba(\Oc_{\pi},\rho)=\infty$, for any $\pi$ in $\an$ of odd
order, except for $\pi=(1 \, 2 \, 3)$ or $\pi=(1 \, 3 \, 2)$ in
$\ac$ -- see Theorem \ref{lemaan}. This last case is particularly
interesting. It corresponds to a ``tetrahedron" rack with constant
cocycle $\omega \in \G_3-1$. The technique in the present paper
does not provide information on the corresponding Nichols algebra.
We also give partial results on pointed Hopf algebras with groups
$\ac$ and $\A_6$, and on Nichols algebras $\toba(\Oc_{\pi},\rho)$,
with $|\pi|$ even. In Section \ref{nichols-dn}, we apply the
technique to conjugacy classes in dihedral groups. It turns out
that it is possible to decide when the associated Nichols algebra
is finite-dimensional in all cases except for $M(\Oc_x, \sgn)$ (if
$n$ is odd), or $M(\Oc_{x}, \sgn \otimes \sgn)$ or $M(\Oc_{x},
\sgn\otimes \varepsilon)$ or $M(\Oc_{xy}, \sgn \otimes \sgn)$ or
$M(\Oc_{xy}, \sgn\otimes \varepsilon)$ (if $n$ is even). See below
for undefined notations. We finally observe in Section
\ref{nichols-other} that there is no finite-dimensional Hopf
algebra with coradical isomorphic to the Hopf algebra $(\ku
\mathbb A_5)^J$ discovered in \cite{Ni}, except for $(\ku \mathbb
A_5)^J$ itself.

\section{Generalities}\label{conventions}

For $s\in G$ we denote by $G^s$ the centralizer of $s$ in $G$. If
$H$ is a subgroup of $G$ and $s\in H$ we will denote $\Oc_{s}^{H}$
the conjugacy class of $s$ in $H$. Sometimes we will write in rack
notations $x\trid y=xyx^{-1}$, $x$, $y\in G$. Also, if $(V,c)$ is
a braided vector space, that is $c\in GL(V\otimes V)$ is a
solution of the braid equation, then $\toba(V)$ denotes its
Nichols algebra.

We denote by $\G_n$ the group of $n$-th roots of 1 in $\ku$.

\subsection{Preliminaries} Let $G$ be a finite group, $\Oc$ a conjugacy class of
$G$, $s\in \Oc$ fixed, $\rho$ an irreducible representation of
$G^s$, $M(\Oc, \rho)$ the corresponding irreducible
Yetter-Drinfeld module. Let $t_1 = s$, \dots, $t_{M}$ be a
numeration of $\Oc$ and let $g_i\in G$ such that $g_i \trid s =
t_i$ for all $1\le i \le M$. Then $M(\Oc, \rho) = \oplus_{1\le i
\le M}g_i\otimes V$. Let $g_iv := g_i\otimes v \in M(\Oc,\rho)$,
$1\le i \le M$, $v\in V$. If $v\in V$ and $1\le i \le M$, then the
coaction and the action of $g\in G$ are given by
$$
\delta(g_iv) = t_i\otimes g_iv, \qquad g\cdot (g_iv) =
g_j(\gamma\cdot v),
$$
where $gg_i = g_j\gamma$, for some $1\le j \le M$ and $\gamma\in
G^s$. The Yetter-Drinfeld module $M(\Oc, \rho)$ is a braided
vector space  with braiding given by
\begin{equation} \label{yd-braiding}
c(g_iv\otimes g_jw) = t_i\cdot(g_jw)\otimes g_iv = g_h(\gamma\cdot
v)\otimes g_iv\end{equation} for any $1\le i,j\le M$, $v,w\in V$,
where $t_ig_j = g_h\gamma$ for unique $h$, $1\le h \le M$ and
$\gamma \in G^s$. Since $s\in Z(G^s)$, the Schur Lemma says that
\begin{equation}\label{schur}
s \text{ acts by a scalar $q_{ss}$ on } V.
\end{equation}

Let $G$ be a finite non-abelian group. Let $\Oc$ be a conjugacy
class of $G$ and let $\rho$ be an irreducible representation of
the centralizer $G^s$ of a fixed $s\in \Oc$. Let $M(\Oc, \rho)$ be
the irreducible Yetter-Drinfeld module corresponding to $(\Oc,
\rho)$ and let $\toba(\Oc, \rho)$ be its Nichols algebra. As
explained in \cite{AZ, AF, G1}, we look for a braided subspace $U$
of $M(\Oc, \rho)$ of diagonal type such that the dimension of the
Nichols algebra $\toba(U)$ is infinite. This implies that the
dimension of $\toba(\Oc, \rho)$ is infinite too.

\begin{lema}\label{trivialbraiding}
If $W$ is a subspace of $V$ such that $c(W\otimes W) = W\otimes W$
and $\dim \toba(W) =\infty$, then $\dim \toba(V) =\infty$.\qed
\end{lema}

\bigbreak Recall that a braided vector space $(V,c)$ is of
\emph{diagonal type} if there exists a basis $v_1, \dots,
v_{\theta}$ of $V$ and non-zero scalars $q_{ij}$, $1\le i,j\le
\theta$, such that $c(v_i\otimes v_j) = q_{ij} v_j\otimes v_i$,
for all $1\le i,j\le \theta$. A braided vector space  $(V,c)$ is
of \emph{Cartan type} if it is of diagonal type and there exists
$a_{ij} \in \Z$, $-|q_{ii}| < a_{ij} \leq 0$ such that
$q_{ij}q_{ji} = q_{ii}^{a_{ij}}$ for all $1\le i\neq j\le \theta$;
by $|q_{ii}|$ we mean $\infty$ if $q_{ii}$ is not a root of 1,
otherwise it means the order of $q_{ii}$ in the multiplicative
group of the units in $\ku$. Set $a_{ii}=2$ for all $1\le i\le
\theta$. Then $(a_{ij})_{1\le i,j\le \theta}$ is a generalized
Cartan matrix.

\begin{theorem}\label{cartantype} (\cite[Th. 4]{H1}, see also
\cite[Th. 1.1]{AS1}). Let $(V,c)$ be a braided vector space of
Cartan type. Then $\dim \toba(V) < \infty$ if and only if the
Cartan matrix is of finite type. \qed
\end{theorem}

We say that $s\in G$ is \emph{real} if it is conjugate to
$s^{-1}$; if $s$ is real, then the conjugacy class of $s$ is also
said to be \emph{real}. We say that $G$ is \emph{real} if any
$s\in G$ is real.

The next application of Theorem \ref{cartantype} was given in
\cite{AZ}. Let $G$ be a finite group, $s\in G$, $\Oc$  the
conjugacy class of $s$, $\rho: G^s \to GL(V)$ irreducible;
$q_{ss}\in \ku^{\times}$ was defined in \eqref{schur}.

\begin{lema}\label{odd}  Assume that $s$ is real. If
$\dim\toba(\Oc, \rho)< \infty$ then $q_{ss} = -1$ and $s$ has even
order.\qed
\end{lema}

If $s^{-1}\neq s$, this is \cite[Lemma 2.2]{AZ}; if $s^2 = \id$
then $q_{ss} = \pm 1$ but $q_{ss} = 1$ is excluded by Lemma
\ref{trivialbraiding}.

The class of real groups includes finite Coxeter groups. Indeed,
all the characters of a finite Coxeter group are real valued, see
subsection \ref{absreal} below, and \cite{BG} for $H_4$.
Therefore, we have:

\begin{theorem}\label{nichols-weyl-impar} Let $G$ be a finite Coxeter group. If $s \in G$ has
odd order, then $\dim\toba(\Oc_s, \rho) = \infty$, for any $\rho
\in \widehat{G^{s}}$. \qed
\end{theorem}

\subsection{Absolutely real groups}\label{absreal}
Let $G$ be a finite group. We say that $s \in G$ is
\emph{absolutely real} if there exists an \emph{involution}
$\sigma$ in $G$ such that $\sigma s \sigma=s^{-1}$. If this
happens, any element in the conjugacy class of $s$ is absolutely
real and we will say that the conjugacy class of $s$ is
\emph{absolutely real}. We say that $G$ is absolutely real if any
$s\in G$ is so. The finite Coxeter groups are absolutely real.
Indeed,
\begin{itemize}
\item[(i)] the dihedral groups are absolutely real, by straightforward
computations.
\item[(ii)] the Weyl groups of semisimple finite dimensional
Lie algebras are absolutely real, by \cite[Th. C (iii), p. 45]{C}.
\item[(iii)] $H_3$ is absolutely real, by Proposition \ref{a5-h4-absreal} below.
\item[(iv)] $H_4$ is absolutely real, we have checked it using GAP3, \cite{Sch97}.
\end{itemize}

\begin{obs}\label{sorites-invreal} Let $G$, $H$ be finite groups. We
note:
\begin{itemize}
    \item $(s,t) \in G\times H$ is absolutely real iff both $s\in G$ and $t\in
    H$ are absolutely real.
    \item $G\times H$ is absolutely real iff both $G$ and $H$ are absolutely real.
    \item Assume $H$ abelian. Then $H$ is absolutely real iff $H$
    has exponent 2, i. e. $H\simeq \Z_2^n$ for some integer $n$.
    \item If $G$ is absolutely real and $H$ is abelian of exponent 2 then
    $G\times H$ is absolutely real.
\end{itemize}
\end{obs}

\bigbreak

We first discuss when an element of $\an$ is absolutely real.
Assume that $\trasp\in\sn$ is of type $(1^{m_1}, 2^{m_2}, \dots,
n^{m_n})$. Then $\trasp\in\an$ iff $\displaystyle\sum_{j \, \,
\text{even}} m_j$ is even.

\begin{lema}\label{an-invreal}
(a). If $m_1 \geq 2$, then $\pi$ is absolutely real in $\an$.

\noindent (b). If $\displaystyle\sum_{h\in \N} (m_{4h} + m_{4h
+3}) $ is even then $\pi$ is absolutely real in $\an$.
\end{lema}

\pf Let $\tau_j := (1\,2\, \dots \, j)$ for some $j$ and take
$$g_j =
\begin{cases} (1 \,\,\, j-1) (2\,\,\, j-2) \cdots (k-1\,\,\, k+1), &\text{ if } j = 2k \text{ is even}, \\
 (1 \,\,\, j-1) (2 \,\,\, j-2) \cdots (k \,\,\, k+1), &\text{ if } j = 2k + 1 \text{ is odd.}
\end{cases}
$$
It is easy to see that $g_j \tau_j g_j = \tau_j^{-1}$, $g_j^2=
\id$ and
$$
\sgn(g_j) =
\begin{cases} (-1)^{k-1}, &\text{ if } j = 2k \text{ is even}, \\
(-1)^{k}, &\text{ if } j = 2k + 1 \text{ is odd.}
\end{cases}
$$

To prove (b), observe that there exists an involution $\sigma\in
\sn$ such that $\sigma \pi\sigma = \pi^{-1}$, which is a product
of ``translations" of the $g_j$'s. Since the sign of $\sigma$ is
$(-1)^{\sum_{h\in \N} (m_{4h} + m_{4h +3})}$, $\sigma \in \an$ iff
$\sum_{h\in \N} (m_{4h} + m_{4h +3})$ is even; (b) follows. We
prove (a). By assumption there are at least two points fixed by
$\pi$, say $n-1$, $n$. By the preceding there exists an involution
$\sigma\in \mathbb S_{n-2}$ such that $\sigma \pi\sigma =
\pi^{-1}$. If $\sigma\in \mathbb A_{n-2} \subset \an$ we are done,
otherwise take $\widetilde \sigma = \sigma (n-1\, n)\in \an$;
$\widetilde \sigma$ is an involution and $\widetilde \sigma
\pi\widetilde \sigma = \pi^{-1}$. \epf

\begin{prop}\label{a5-h4-absreal}
The groups $\aco$ and $H_3$ are absolutely real.
\end{prop}

\pf The type of $\pi\in \aco$ is either $(1^5)$, $(1^2,3^1)$,
$(1,2^2)$ or $(5^1)$; in the first two cases $\pi$ is absolutely
real by Lemma \ref{an-invreal} part (a), in the last two by part
(b). Since $H_3\simeq \aco \times \Z_2$ (see \cite[Section
2.13]{Hu}), then the Coxeter group is absolutely real by Remark
\ref{sorites-invreal}. \epf

\subsection{Generalizations of Lemma \ref{odd}}
The next two Lemmata are variations of \cite[Lemma 2.2]{AZ}. A
result in the same spirit appears in \cite{FGV}. We deal with
elements $s$ having a power in $\Oc$, the conjugacy class of $s$.
Clearly, if $s^j=\sigma s \sigma^{-1}$ is in $\Oc$, then
$s^{j^l}=\sigma^l s \sigma^{-l}$ is in $\Oc$, for every $l$. So,
$s^{j^{|\sigma|}}=s$; this implies that $|s|$ divides
$j^{|\sigma|}-1$. Hence
\begin{align}\label{ecN}
N  \quad \text{divides} \quad j^{|\sigma|}-1,
\end{align}
with $N:=|q_{ss}|$, recall \eqref{schur}.

\begin{lema}\label{lemagen2}
Let $G$ be a finite group, $s\in G$, $\Oc$ the conjugacy class of
$s$ and $\rho \in \widehat{G^s}$. Assume that there exists an
integer $j$ such that $s$, $s^j$ and $s^{j^2}$ are distinct
elements and $s^j$ is in $\Oc$. If $\dim\toba(\Oc,\rho)<\infty$,
then $s$ has even order and $q_{ss}=-1$.
\end{lema}

\pf We assume that $\dim\toba(\Oc,\rho)<\infty$, thus $N>1$. It is
easy to see that
\begin{align}\label{ecrel}
\sigma^{-h} s^{j^l} \sigma^h=(\sigma^{-h} s
\sigma^h)^{j^l}=(s^{j^{|\sigma|-h}})^{j^l}=s^{j^{|\sigma|-h+l}},
\end{align}
for every $l$, $h$. We will call $t_l:=s^{j^l}$,
$g_l:=\sigma^{l}$, $l=0$, $1$, $2$; so $t_l=g_l s g_l^{-1}$, for
$l=0$, $1$, $2$. The other relations between $t_l$'s and $g_h$'s
are obtained from \eqref{ecrel}. For $v \in V-0$ and $l=1$ or $2$,
we define $W_l:=\ku-\text{span of } \{g_0v,g_lv \}$. Hence, $W_l$
is a braided vector subspace of $M(\Oc,\rho)$ of Cartan type with
$$\mathcal Q_l =\begin{pmatrix}
q_{ss} & q_{ss}^{j^{|\sigma|-l}} \\ q_{ss}^{j^l} & q_{ss}
\end{pmatrix},\qquad
\mathcal A_l =\begin{pmatrix} 2 & a_{12}(l) \\ a_{21}(l) & 2
\end{pmatrix},
$$
where $a_{12}(l)=a_{21}(l)\equiv j^{|\sigma|-l}+j^l \mod(N)$.
Since $\dim\toba(\Oc,\rho)<\infty$, we
have that $a_{12}(l)=a_{21}(l)= 0$ or $-1$. We consider now two cases.\\
(i) Let us suppose that $a_{12}(1)=a_{21}(1)= 0$. This implies
that $j^{|\sigma|-1}+j \equiv 0 \mod(N)$. Since $N$ divides
$j^{|\sigma|}-1$, we have that $N$ divides $j^2+1$. We consider
now two possibilities.
\begin{itemize}
\item Assume that $a_{12}(2)=a_{21}(2)= 0$. Then $j^{|\sigma|-2}+j^2 \equiv 0
\mod(N)$. Since $N$ divides $j^{|\sigma|}-1$, we have that $N$
divides $j^4+1$. So, $-1 \equiv 1 \mod (N)$; hence the result
follows.
\item Assume that $a_{12}(2)=a_{21}(2)= -1$. Then $j^{|\sigma|-2}+j^2 \equiv -1
\mod(N)$. We can see that $N$ divides $j^4+j^2+1$. This implies
that $N$ divides $1$, a contradiction.
\end{itemize}
(ii) Let us suppose that $a_{12}(1)=a_{21}(1)= -1$. This implies
that $j^{|\sigma|-1}+j \equiv -1 \mod(N)$. Since $N$ divides
$j^{|\sigma|}-1$, we have that $N$ divides $j^2+j+1$. We consider
now two possibilities.
\begin{itemize}
\item Assume that $a_{12}(2)=a_{21}(2)= 0$. Then $j^{|\sigma|-2}+j^2 \equiv 0
\mod(N)$. So, $N$ divides $j^4+1$. It is easy to see that $N$
divides $j^2$. Since $j$ and $|s|$ are relatively prime, $N$ must
be $1$, a contradiction.
\item Assume that $a_{12}(2)=a_{21}(2)= -1$. This means that the subspace
$\widetilde{W}:=\ku-\text{span of } \{g_0v,g_1v,g_2v \}$ of
$M(\Oc,\rho)$ is of Cartan type with
$$\mathcal Q =\begin{pmatrix}
q_{ss} & q_{ss}^{j^{|\sigma|-1}} & q_{ss}^{j^{|\sigma|-2}}\\
q_{ss}^j & q_{ss} & q_{ss}^{j^{|\sigma|-1}}\\
q_{ss}^{j^2} & q_{ss}^j & q_{ss}
\end{pmatrix}, \qquad
\mathcal A =\begin{pmatrix} 2 & -1 & -1 \\ -1 & 2 &-1 \\ -1 & -1
&2
\end{pmatrix}.
$$
By Theorem \ref{cartantype}, we have that
$\dim\toba(\Oc,\rho)=\infty$, a contradiction.
\end{itemize}
This concludes the proof. \epf

\begin{lema}\label{lemagen1}
Let $G$ be a finite group, $s\in G$, $\Oc$ the conjugacy class of
$s$ and $\rho=(\rho,V) \in \widehat{G^s}$ such that
$\dim\toba(\Oc,\rho)<\infty$. Assume that there exists an integer
$j$ such that $s^j\neq s$ and $s^j$ is in $\Oc$. \begin{itemize}
\item[(a)] If $\deg \rho >1$, then $s$ has even order and
$q_{ss}=-1$.
\item[(b)] If $\deg \rho =1$, then either $s$ has even order and $q_{ss}=-1$, or $q_{ss}$
$\in \G_3-1$.
\end{itemize}
\end{lema}

\pf We will proceed and use the notation as in the proof of Lemma
\ref{lemagen2}. If $s^{j^2}\neq s$, then the result follows by
Lemma \ref{lemagen2}. Assume that $s^{j^2}= s$. This implies that
$|s|$ divides $j^2-1$, so $N$ divides $j^2-1$.

(a) Let $v_1$ and $v_2$ in $V$ linearly independent and let
$W=\ku$-- span of $\{g_0 v_1, \,  g_0 v_2, \,  g_1 v_1, \,  g_1
v_2 \}$, with $g_0:=\id$ and $g_1:=\sigma$. Thus $W$ is a braided
vector subspace of $M(\Oc,\rho)$ of Cartan type with
$$\mathcal Q =\begin{pmatrix}
q_{ss} & q_{ss} & q_{ss}^{j^{|\sigma|-1}} & q_{ss}^{j^{|\sigma|-1}}\\
q_{ss} & q_{ss} & q_{ss}^{j^{|\sigma|-1}} & q_{ss}^{j^{|\sigma|-1}}\\
q_{ss}^j  & q_{ss}^j & q_{ss} & q_{ss}\\
q_{ss}^j  & q_{ss}^j & q_{ss} & q_{ss}
\end{pmatrix}, \qquad
\mathcal A =\begin{pmatrix} 2 & a_{12} & a_{13} & a_{14} \\
a_{21} & 2 & a_{23}  & a_{24} \\
a_{31} & a_{32} & 2  & a_{34} \\
a_{41} & a_{42} & a_{43}  & 2
\end{pmatrix},
$$
where $a_{ij}=a_{ji}$, $i\neq j$, $a_{12}\equiv 2 \equiv a_{34}
\mod (N)$, $a_{13}=a_{14}=a_{23}=a_{24}$ and
$$
a_{13} \equiv j^{|\sigma|-1}+j \mod(N).
$$
If $a_{12}=0$ or $a_{34}=0$, then $N$ divides $2$ and the result
follows. Besides, if $a_{13}=0$, then $j^{|\sigma|-1}+j\equiv 0
\mod(N)$; this implies that $N$ divides $j^2+1$, so $N$ divides 2
and the result follows. On the other hand, if $a_{ij}=-1$, for all
$i,j$, we have that the matrix $\mathcal A$ is not of finite type;
hence $\dim \toba (\Oc, \rho)=\infty$, from Theorem
\ref{cartantype}. This is a contradiction by hypothesis.
Therefore, (a) is proved.

(b) For $v \in V-0$ we define $W:=\ku-\text{span of } \{g_0v,g_1v
\}$, with $g_0:=\id$ and $g_1:=\sigma$. Hence, $W$ is a braided
vector subspace of $M(\Oc,\rho)$ of Cartan type with
$$\mathcal Q =\begin{pmatrix}
q_{ss} & q_{ss}^{j^{|\sigma|-1}} \\ q_{ss}^{j} & q_{ss}
\end{pmatrix}, \quad
\mathcal A =\begin{pmatrix} 2 & a_{12} \\ a_{21} & 2
\end{pmatrix},
$$
where $a_{12}\equiv j^{|\sigma|-1}+j \mod(N)$. Since
$\dim\toba(\Oc,\rho)<\infty$, we
have that $a_{12}= 0$ or $-1$. We consider now two possibilities.\\
(i) Assume that $a_{12}= 0$. This implies that $j^{|\sigma|-1}+j
\equiv 0 \mod(N)$. Since $N$ divides $j^{|\sigma|}-1$, we have
that $N$ divides $j^2+1$. Thus, $N$ divides $2$; hence $N=2$, and the result follows.\\
(ii) Assume that $a_{12}= -1$. This implies that $j^{|\sigma|-1}+j
\equiv -1 \mod(N)$. Since $N$ divides $j^{|\sigma|}-1$, we have
that $N$ divides $j^2+j+1$. Hence, $N$ divides $j+2$. If $p$ is a
prime divisor of $N$, then $p$ divides $j-1$ or $j+1$, because $N$
divides $j^2-1$. If $p$ divides $j+1$, then $p$ divides $1$, a
contradiction. So, $N$ divides $j-1$. Hence, $N$ divides $3$, i.e.
$N=3$ and the result follows. \epf

\section{On Nichols algebras over $\mathbb A_n$}
\label{nichols-an}

We recall that we will denote $\Oc$ or $\Oc_{\pi}$ the conjugacy
class of an element $\pi$ in $\an$, and $\rho$ in
$\widehat{\an^{\pi}}$, a representative of an isomorphism class of
irreducible representations of $\an^{\pi}$. We want to determine
pairs $(\Oc,\rho)$, for which $\dim \toba(\Oc,\rho)=\infty$,
following the strategy given in \cite{AZ,AF}; see also \cite{G1}.

The following is a helpful criterion to decide when a conjugacy
class of an even permutation $\pi$ in $\sn$ splits in $\mathbb
A_n$.

\begin{prop}\label{split}\cite[Proposition 12.17]{JL}
Let $\pi\in \an$, with $n>1$.
\begin{itemize}
\item[(1)] If $\pi$ commutes with some odd permutation in $\sn$, then
$\Oc_{\pi}^{\an}=\Oc_{\pi}^{\sn}$ and $[\sn^{\pi}:\an^{\pi}]=2$.
\item[(2)] If $\pi$ does not commute with any odd permutation in
$\sn$, then $\Oc_{\pi}^{\sn}$ splits into two conjugacy classes in
$\mathbb A_n$ of equal size, with representatives $\pi$ and $(1 \,
2)\pi(1 \, 2)$, and $\sn^{\pi}=\an^{\pi}$.\qed
\end{itemize}
\end{prop}

\begin{rmk}\label{obser}
(i) Notice that if $\pi$ satisfies (1) of Proposition \ref{split},
then $\pi$ is real. The reciprocal is not true, e.g. consider
$\tau_5=(1 \, 2\, 3\, 4\, 5)$ in $\aco$.

(ii) One can see that if $\pi$ in $\an$ is of type
$(1^{m_1},2^{m_2},\dots,n^{m_n})$, then $\pi$ satisfies (2) of
Proposition \ref{split} if and only if $m_1=0$ or $1$, $m_{2h}=0$
and $m_{2h+1}\leq 1$, for all $h\geq 1$. Thus, if $\pi \in \an$
has even order, then $\pi$ is real.
\end{rmk}

We state the main Theorem of the section.

\begin{theorem}\label{lemaan}
Let $\pi\in \an$ and $\rho \in \widehat{\an^{\pi}}$. Assume that
$\pi$ is neither $(1 \, 2 \, 3)$ nor $(1 \, 3 \, 2)$ in $\ac$. If
$\dim\toba(\Oc_{\pi},\rho) < \infty$, then $\pi$ has even order
and $q_{\pi \pi}=-1$.
\end{theorem}
\pf If $|\pi|$ is even the result follows by Lemma \ref{odd} and
Remark \ref{obser} (ii). Let us suppose that $|\pi|\geq 5$ and odd
. If $\pi^{-1}$ is in $\Oc_{\pi}$, then the result follows by
Lemma \ref{odd}. Assume that $\pi^{-1}\not \in \Oc_{\pi}$. We
consider two cases.

(i) If $\pi^2 \in \Oc_{\pi}$, then $\pi^4$ is in $\Oc_{\pi}$, and
$\pi^4\neq \pi^2$ because $|\pi| \geq 5$. Hence, the result
follows from Lemma \ref{lemagen2}.

(ii) Assume that $\pi^2\not \in\Oc_{\pi}$. We know that there
exist $\sigma$ and $\sigma'$ in $\sn$, necessarily odd
permutations, such that $\pi^{-1}=\sigma \pi \sigma^{-1}$ and
$\pi^{2}=\sigma' \pi \sigma'^{-1}$. Then $\sigma''=\sigma\sigma'
\in \an$ and $\pi^{-2}=\sigma'' \pi \sigma''^{-1}$; so $\pi^{-2}$
is in $\Oc_{\pi}$. This implies that $\pi^4$ is in $\Oc_{\pi}$,
and $\pi^4\neq \pi^{-2}$ because $5 \leq |\pi|$ is odd. Now, the
result follows from Lemma \ref{lemagen2}.

Finally, let us suppose that $|\pi|=3$, with type $(1^a,3^b)$. If
$a\geq 2$ or $b\geq 2$, then $\pi$ is real, by Lemma
\ref{an-invreal} (a) and Remark \ref{obser}, respectively. Hence,
the result follows by Lemma \ref{odd}. This concludes the proof.
\epf

\subsection{Case $\at$}
Obviously, $\at \simeq \Z_3$; thus $\at$ is not real. This case
was considered in \cite[Theorem 1.3]{AS1}.

\subsection{Case $\ac$}

It is straightforward to check that $\ac$ is not real, since $(1
\, 2 \, 3)$ is not real in $\ac$. Let $\pi$ in $\ac$; then the
type of $\pi$ may be $(1^4)$, $(2^2)$ or $(1,3)$. If the type of
$\pi$ is $(1^4)$, then $\dim \toba(\Oc_{\pi}, \rho)=\infty$, for
any $\rho$ in $\widehat{\ac}$, by Lemma \ref{trivialbraiding}. If
the type of $\pi$ is $(1,3)$, then $\pi$ is not real; moreover we
have
$$\Oc_{(1 \, \, 2 \,\, 3)}=\{(1 \, \, 2 \,\, 3),(1 \, \, 3 \, \,
4), (1 \, \, 4 \,\, 2),(2 \, \, 4 \,\, 3)\},$$
$$\Oc_{(1 \, \, 3 \,\, 2)}=\{(1 \, \, 3 \,\, 2),(1 \, \, 2 \,\,
4), (1 \, \, 4 \,\, 3),(2 \, \, 3 \,\, 4)\},$$ and
$\ac^{\pi}=\langle \pi \rangle\simeq \Z_3$. If $\rho \in
\widehat{\ac^{\pi}}$ is trivial, then $\dim \toba(\Oc_{\pi},
\rho)=\infty$; otherwise it is not known.

The following result is a variation of \cite[Theorem 2.7]{AZ}.

\begin{prop}\label{a4}
Let $\pi$ in $\ac$ of type $(2^2)$. Then $\dim \toba(\Oc_{\pi},
\rho)=\infty$, for every $\rho$ in $\widehat{\ac^{\pi}}$.
\end{prop}

\pf We can assume that $\pi=(1 \, 2)(3 \, 4)$. If we call
$t_1:=\pi$, $t_2:=(1 \, 3)(2 \, 4)$ and $t_3:=(1 \, 4) (2 \, 3)$,
then $\Oc_{\pi}= \{ t_1 , t_2 , t_3 \}$ and $\ac^{\pi}=\langle t_1
\rangle \times \langle t_2 \rangle \simeq \Z_2 \times  \Z_2$. If
$g_1=\id$, $g_2=(1 \, 3 \, 2)$ and $g_3=(1 \, 2 \, 3)$, then
$t_j=g_j \pi g_j^{-1}$, $j=1,2,3$, and
\begin{align*}
t_1g_2&=g_2t_3,  \quad &t_2g_1&=g_1t_2, \quad &t_3g_1&=g_1t_3,\\
t_1g_3&=g_3t_2,  \quad &t_2g_3&=g_3t_3,  \quad &t_3g_2&=g_2t_2.
\end{align*}

Let $\rho$ in $\widehat{\ac^{\pi}}$ and $M(\Oc_{\pi},\rho):=g_1v
\oplus g_2v \oplus g_3v$, where $\langle v \rangle$ is the vector
space affording $\rho$. Thus $M(\Oc_{\pi},\rho)$ is a braided
vector space with braiding given by -- see \eqref{yd-braiding}--
$c(g_j v \otimes g_j v)= g_j t_1 \cdot v \otimes g_j v$ and $c(g_j
v \otimes g_1 v)= g_1 t_j \cdot v \otimes g_j v$, $j=1,2,3$ and
\begin{align*}
c(g_1 v \otimes g_2v)&=g_2 t_3 \cdot v \otimes g_1v, \qquad
&c(g_1v \otimes g_3v)&=g_3 t_2 \cdot v \otimes g_1v, \\
c(g_2 v \otimes g_3v)&=g_3 t_3 \cdot v \otimes g_2v, \qquad
&c(g_3v \otimes g_2v)&=g_2 t_2 \cdot v \otimes g_3v.
\end{align*}

Clearly, $\dim\toba(\Oc_{\pi},\varepsilon \otimes
\varepsilon)=\dim\toba(\Oc_{\pi},\varepsilon \otimes
\sgn)=\infty$, by Lemma \ref{trivialbraiding}. If we consider
$\rho=\sgn \otimes \varepsilon$ (resp. $\sgn \otimes \sgn$), then
$M(\Oc_{\pi},\rho)$ is of Cartan type with matrix of coefficients
$(q_{ij})_{ij}$ given by
$$\mathcal Q =\begin{pmatrix}
-1 & -1 & 1\\ 1 & -1 & -1\\
-1 & 1 & -1
\end{pmatrix},\qquad
(\text{ resp. } \mathcal Q =\begin{pmatrix}
-1 & 1 & -1\\ -1 & -1 & 1\\
1 & -1 & -1
\end{pmatrix}).
$$
In both cases the Cartan matrix is $\mathcal A =\begin{pmatrix}
2 & -1 & -1\\ -1 & 2 & -1\\
-1 & -1 & 2
\end{pmatrix}$.
Therefore, $\dim\toba(\Oc_{\pi}, \rho)=\infty$, by Theorem
\ref{cartantype}. \epf

\medbreak
\subsection{Case $\aco$}
Here is the key step in the consideration of this case.

\begin{lema}\label{a5}
Let $\pi \in \aco$. Then $\dim \toba(\Oc_{\pi}, \rho)=\infty$, for
every $\rho$ in $\widehat{\aco^{\pi}}$.
\end{lema}

\pf Let $\pi \in \aco$. If the type of $\pi$ is either $(1^5)$,
$(1^2,3)$ or $(5)$, we have that $\dim \toba(\Oc_{\pi},
\rho)=\infty$, by Lemma \ref{odd} and Proposition
\ref{a5-h4-absreal}. Let us assume that the type of $\pi$ is
$(2^2)$. For $j=1,2,3$, let $t_j$ and $g_j$ be as in the proof of
Proposition \ref{a4}. By Proposition \ref{split} and
straightforward computations, we have that
$\Oc_{\pi}^{\aco}=\Oc_{\pi}^{\mathbb S_5}$ and $\aco^{\pi}=\langle
t_1 \rangle \times \langle t_2 \rangle \simeq \Z_2 \times  \Z_2$.
Notice that $t_j \in \Oc_{\pi}^{\aco}$, $j=1,2,3$. Let $\rho\in
\widehat{\aco^{\pi}}$ and $W:=g_1v \oplus g_2v \oplus g_3v$, where
$\langle v \rangle$ is the vector space affording $\rho$; then $W$
is a braided vector subspace of $M(\Oc_{\pi}, \rho)$. Therefore,
$\dim \toba(\Oc_{\pi}, \rho)=\infty$, by the same argument as in
the proof of Proposition \ref{a4}. \epf

As an immediate consequence of Lemma \ref{a5} we have the
following result.

\begin{theorem}\label{mainteor}
Any finite-dimensional complex pointed Hopf algebra $H$ with
$G(H)\simeq \aco$ is necessarily isomorphic to the group algebra
of $\aco$.
\end{theorem}

\pf Let $H$ be a complex pointed Hopf algebra with $G(H)\simeq
\aco$. Let $M \in {}_{\ku \aco}^{\ku \aco}\mathcal{YD}$ be the
infinitesimal braiding of $H$ -see \cite{AS-cambr}. Assume that
$H\neq \ku \aco$; thus $M\neq 0$. Let $N\subset M$ be an
irreducible submodule. Then $\dim\toba(N)=\infty$, by Lemma
\ref{a5}. Hence, $\dim\toba(M)=\infty$ and $\dim H=\infty$. \epf

\subsection{Case $\A_6$}
Let $\pi$ be in $\A_6$. If the type of $\pi$ is $(1^6)$,
$(1^2,2^2)$, $(1^3,3)$, $(3^2)$ or $(1,5)$, then $\pi$ is
absolutely real by Lemma \ref{an-invreal}, and if the type of
$\pi$ is $(2,4)$, then $\pi$ is real because it has even order --
see Remark \ref{obser} (ii). Hence, $\A_6$ is a real group. We
summarize our results in the following statement.

\begin{theorem}
Let $M(\Oc,\rho)$ be an irreducible Yetter-Drinfeld module over
$\ku \mathbb A_6$, corresponding to a pair $(\Oc,\rho)$. If
$\dim\toba(\Oc,\rho)<\infty$, then $\Oc=\Oc_{\pi}$, with $\pi=(1
\, 2)(3 \, 4 \, 5 \, 6)$, and $\rho=\sgn \in \widehat{\Z_4}$.
\end{theorem}

\begin{obs}
In this Theorem we do not claim that the condition is sufficient.
\end{obs}

\pf Let $\pi$ be in $\A_6$. If the type of $\pi$ is
\begin{itemize}
\item $(1^6)$, then $\dim
\toba(\Oc_{\pi}, \rho)=\infty$, for any $\rho$ in
$\widehat{\A_{6}^{\pi}}$, by Lemma \ref{trivialbraiding}.
\item $(1^3,3)$, $(3^2)$ or $(1,5)$, then $\dim
\toba(\Oc_{\pi}, \rho)=\infty$, for any $\rho$ in
$\widehat{\A_{6}^{\pi}}$, by Lemma \ref{odd}.
\end{itemize}

Let us suppose that the type of $\pi$ is $(1^2,2^2)$; we can
assume that $\pi=(1 \, 2)(3\, 4)$. It is easy to check that
$$\A_6^{\pi}=\langle a:=(3 \, 4)(5 \, 6), \, b:=(1 \, 3 \, 2 \, 4)(5 \, 6)\rangle \simeq \mathbb D_4.$$
Notice that $\pi=b^2$. It is known that $\widehat{\mathbb
D_4}=\{\rho_1, \, \rho_2,
 \, \rho_3, \, \rho_4, \, \rho_5 \}$, where $\rho_j$, $j=1$, $2$, $3$ and
 $4$, are the following characters
\begin{align*}
\rho_1(a)&= 1, \quad &\rho_2(a)&= -1, \quad &\rho_3(a)&= 1, \quad &\rho_4(a)&= -1,\\
\rho_1(b)&= 1, \quad &\rho_2(b)&= 1, \quad &\rho_3(b)&= -1, \quad
&\rho_4(b)&= -1,
\end{align*}
and $\rho_5$ is the $2$-dimensional representation given by
\begin{align*}
\rho_5(a)=\begin{pmatrix} 0 & 1 \\ 1 & 0
\end{pmatrix}, \qquad
\rho_5(b)=\begin{pmatrix} i & 0 \\ 0 & -i
\end{pmatrix}.
\end{align*}
It is clear that $\rho_j(\pi)=1$, $j=1$, $2$, $3$ and $4$. Then
$\dim\toba(\Oc_{\pi},\rho)=\infty$, by Lemma
\ref{trivialbraiding}. Let us consider now that $\rho=\rho_5$. We
define $t_1:=(1 \, 2)(3 \, 4)$, $t_2:=(1 \, 3)(2 \, 4)$, $t_3:=(1
\, 4)(2 \, 3)$, $g_1:=\id$, $g_2:=(1 \, 3 \, 2)$ and $g_3:=(1  \,
2 \, 3)$. It is clear that
\begin{align*}
\rho(t_1)=\begin{pmatrix} -1 & 0 \\ 0 & -1
\end{pmatrix}, \qquad
\rho(t_2)=\begin{pmatrix} 0 & i \\ -i & 0
\end{pmatrix}=-\rho(t_3).
\end{align*}
If $v_1:=\begin{pmatrix} i \\ 1\end{pmatrix}$ we have that
$\rho(t_1)(v_1)=-v_1$ and $\rho(t_2)(v_1)=v_1=-\rho(t_3)(v_1)$. We
define $W:=\ku-\text{span of } \{g_1v_1,g_2v_1,g_3v_1 \}$. Then
$W$ is a braiding subspace of $M(\Oc_{\pi},\rho)$ of Cartan type
with
\begin{align*}
\mathcal Q=\begin{pmatrix} -1 & -1 & 1 \\ 1 & -1 & -1
\\ -1 & 1 & -1
\end{pmatrix}, \qquad
\mathcal A=\begin{pmatrix} 2 & -1 & -1 \\ -1 & 2 & -1
\\ -1 & -1 & 2
\end{pmatrix}.
\end{align*}
Since $\mathcal A$ is not of finite type we have that
$\dim\toba(\Oc_{\pi},\rho)=\infty$, by Theorem \ref{cartantype}.

Finally, let us assume that the type of $\pi$ is $(2,4)$. Then
$\Oc_{\pi}$ has $90$ elements and $\A_{6}^{\pi}=\langle \pi
\rangle \simeq \Z_4$. We call
$\widehat{\Z_4}=\{\chi_0,\chi_1,\chi_2,\chi_3 \}$, where
$\chi_l(\pi)=i^l$, $l=0$, $1$, $2$, $3$. It is clear that if
$\rho=\chi_l$, with $l=0$, $1$ or $3$, then $\rho(\pi)\neq -1$.
This implies that $\dim\toba(\Oc_{\pi},\rho)=\infty$, by Lemma
\ref{odd}. \epf

\begin{obs}
We can see that every maximal abelian subrack of
$\Oc_{(12)(3456)}$ has two elements. Hence,
$M(\Oc_{(12)(3456)},\rho)$ is a negative braided space in the
sense of \cite{AF}.
\end{obs}

\medbreak
\subsection{Case $\am$, $m\geq 7$}

Let $\pi  \in \A_m$, with $|\pi|$ even. We now investigate the
Nichols algebras associated with $\pi$ by reduction to the
analogous study for the orbit of $\pi$ in $\sm$, \cite{AF}. By
Remark \ref{obser} (ii), $\Oc_{\pi}^{\am}=\Oc_{\pi}^{\sm}$ and
$[\sm^{\pi}:\am^{\pi}]=2$. So, we can determinate the irreducible
representations of $\am^{\pi}$ from those of $\sm^{\pi}$. We know
that if the type of $\pi$ is $(1^{b_1},2^{b_2},\dots,m^{b_m})$,
then $\sm^{\pi}= T_1 \cdots T_m, $ with $T_i \simeq
\Z_i^{b_i}\rtimes \mathbb S_{b_i}$, $1\leq i\leq m$.

\medbreak {\bf{Some generalities and notation.}} Let $G$ be a
finite group, $H$ a subgroup of $G$ of index two, and $\eta$ a
representation of $G$. It is easy to see that
$$ \eta '(g):=\begin{cases}\eta(g) &, \text{if $g \in H$},\\
-\eta(g) &,\text{if $g \in G \setminus H$},
\end{cases}
$$
defines a new representation of $G$. Notice that
$\Res_{H}^{G}\eta=\Res_{H}^{G}\eta '$. On the other hand, any
representation $\rho$ of $H$ defines a representation
$\overline{\rho}$ of $H$, call it the \emph{conjugate
representation of $\rho$}, given by
$\overline{\rho}(h):=\rho(ghg^{-1})$, for every $h \in H$, where
$g$ is an arbitrary fixed element in $G\setminus H$. Since $g$ is
unique up to multiplication by an element of $H$, the conjugate
representation is unique up to isomorphism.

Let $s\in H$ such that $\Oc_{s}^{H}=\Oc_{s}^{G}$; thus
$[G^s:H^s]=2$. Let $\eta$ in $\widehat{G^s}$. Then we have two
cases:
\begin{itemize}
\item[(i)] $\eta \not\simeq \eta '$. If
$\rho:=\Res_{H^s}^{G^s} \eta$, then $\rho \in \widehat{H^s}$,
$\rho\simeq \overline{\rho}$ and $\Ind_{H^s}^{G^s} \rho \simeq
\eta \oplus \eta'$.
\item[(ii)] $\eta \simeq \eta '$. We have that
$\Res_{H^s}^{G^s} \eta \simeq \rho \oplus \overline{\rho}$ and
$\Ind_{H^s}^{G^s} \rho \simeq \eta \simeq \Ind_{H^s}^{G^s}
\overline{\rho}$.
\end{itemize}
Moreover, if $\rho$ is an irreducible representation of $H^s$,
then $\rho$ is a restriction of some $\eta \in \widehat{G^s}$ or
is a direct summand of $\Res_{H^s}^{G^s} \eta$ as in (ii), see
\cite[Ch. 5]{FH}.

\begin{obs} If $\eta \in \widehat{G^s}$ and $\rho:=\Res_{H^s}^{G^s} \eta$, it is easy to check
that
\begin{align}
M(\Oc^{G}_{s},\eta)&\simeq M(\Oc^{H}_{s},\rho), & &\text{for the case (i)},\\
M(\Oc^{G}_{s},\eta)&\simeq M(\Oc^{H}_{s},\rho) \oplus
M(\Oc^{H}_{s},\overline{\rho}), & & \text{for the case (ii)},
\end{align}
as braided vector spaces.
\end{obs}

We apply these observations to the case $G=\sm$ and $H=\am$. We
use some notations given in \cite[Section II.D]{AF}.

\begin{lema} Assume that the type of $\pi$ is $((2r)^n)$, with $r\geq 1$
and $n$ even. Let $\rho$ in $\widehat{\A_{m}^{\pi}}$, with
$m=2rn$.

\begin{itemize}
\item[(a)] If $q_{\pi \pi}\neq -1$, then $\dim\toba(\Oc_{\pi}, \rho) = \infty$.
\item[(b)] If $\rho \simeq \overline{\rho}$ and $q_{\pi \pi}= -1$, then
\begin{itemize}
\item[(I)] if $r=1$, then $\dim \toba(\Oc, \rho)= \infty$.
\item[(II)] Assume that $r > 1$. If $\deg \rho>1$, then $\dim \toba(\Oc, \rho)= \infty$.
Assume that $\deg \rho=1$. If $\rho=\chi_{r,\dots,r} \otimes \mu$,
with $r$ even or odd, or if $\rho=\chi_{c,\dots,c} \otimes \mu$,
with $r$ even and $c=\frac{r}{2}$ or $\frac{3r}{2}$, where
$\mu=\varepsilon$ or $\sgn$, then the braiding is negative;
otherwise, $\dim \toba(\Oc, \rho)= \infty$.
\end{itemize}
\end{itemize}
\end{lema}

\pf (a) follows by Remark \ref{obser} (ii) and Lemma \ref{odd}.
(b). Since \emph{$\rho \simeq \overline{\rho}$},
$\rho=\Res^{\sm^{\pi}}_{\am^{\pi}}(\eta)$, with $\eta\in
\widehat{\sm^{\pi}}$, $\eta \not \simeq \eta '$ and
$\Ind^{\sm^{\pi}}_{\am^{\pi}}\rho \simeq \eta \oplus \eta '$.
Notice that $\eta(\pi)=-\Id$ because $\rho(\pi)=-\Id$. Now, as the
racks are the same, i.e. $\Oc_{\pi}^{\am}=\Oc_{\pi}^{\sm}$, we can
apply \cite[Theorem 1]{AF}. \epf

\begin{obs} Keep the notation of the Lemma. If $\rho$ is not isomorphic to its conjugate
representation $\overline{\rho}$, then there exists $\eta \in
\widehat{\sm^{\pi}}$ such that
$\Res^{\sm^{\pi}}_{\am^{\pi}}(\eta)=\rho \oplus \overline{\rho}$
and $\Ind^{\sm^{\pi}}_{\am^{\pi}}\rho \simeq \eta \simeq \eta '
\simeq \Ind^{\sm^{\pi}}_{\am^{\pi}}\overline{\rho}$. Clearly,
$\eta(\pi)$ and $\overline{\rho}(\pi)$ act by scalar $-1$, and we
have that $ M(\Oc^{\sm}_{\pi},\eta)\simeq  M(\Oc^{\am}_{\pi},\rho)
\oplus M(\Oc^{\am}_{\pi},\overline{\rho})$ as braided vector
spaces. We do not get new information with the techniques
available today.

\end{obs}

\bigbreak
\section{On Nichols algebras over $\dn$}\label{nichols-dn}

We fix the notation: the dihedral group $\dn$ of order $2n$ is
generated by $x$ and $y$ with defining relations $x^2 = e = y^n$
and $xyx = y^{-1}$. Let $\omega$ be a primitive $n$-th root of 1
and let $\chi$ be the character of $\langle y \rangle$,
$\chi(y)=\omega$. If $s \in \dn$ then we denote the conjugacy
class by $\Oc_s^n$ or simply $\Oc_s$.

\begin{theorem}\label{dn} Let $M(\Oc, \rho)$ be the
irreducible Yetter-Drinfeld module over $\ku\dn$ corresponding to
a pair $(\Oc, \rho)$. Assume that its Nichols algebra $\toba(\Oc,
\rho)$ is finite-dimensional.

\begin{enumerate}
\item[(a)] If $n$ is odd, then $(\Oc, \rho) = (\Oc_{x},
\sgn)$, where $\sgn \in \widehat{\dn^x}$, $\dn^x=\langle x \rangle
\simeq \mathbb Z_2$.

\item[(b)] If $n =2m$ is even, then $(\Oc, \rho)$ is one of the following:

\begin{itemize}
\item[(i)] $(\Oc_{y^{m}}, \rho)$ where $\rho\in \widehat{\dn}$
satisfies $\rho(y^{m}) = -1$.

\item[(ii)] $(\Oc_{y^{h}}, \chi^j)$ where $1\leq h\leq n-1$, $h\neq m$
and $\omega^{hj} = -1$.

\item[(iii)] $(\Oc_{x}, \sgn \otimes \sgn)$ or $(\Oc_{x}, \sgn\otimes
\varepsilon)$, where $\sgn \otimes \sgn$, $\sgn\otimes \varepsilon
\in \widehat{\dn^x}$, $\dn^x=\langle x \rangle \oplus \langle y^m
\rangle \simeq \mathbb Z_2 \times \mathbb Z_2$.

\item[(iv)] $(\Oc_{xy}, \sgn \otimes \sgn)$ or $(\Oc_{xy}, \sgn\otimes
\varepsilon)$, where $\sgn \otimes \sgn$, $\sgn\otimes \varepsilon
\in \widehat{\dn^{xy}}$, $\dn^{xy}=\langle xy \rangle \oplus
\langle y^m \rangle \simeq \mathbb Z_2 \times \mathbb Z_2$.

\end{itemize}

\end{enumerate}

In the cases (i) and (ii) the dimension is finite. In the cases
(iii) and (iv), the braiding is negative in the sense of
\cite{AF}.
\end{theorem}

\begin{obs}\label{obs:isoracks} There are isomorphisms of braided vector spaces
\begin{align*}M(\Oc_{x}, \sgn \otimes \sgn) &\simeq M(\Oc_{xy}, \sgn
\otimes \sgn), \\ M(\Oc_{x}, \sgn\otimes \varepsilon) &\simeq
M(\Oc_{xy}, \sgn \otimes \varepsilon).
\end{align*}\end{obs}

\begin{obs}\label{cociente} Assume for simplicity that $n$ is odd and that $n= de$,
where $d$, $e$ are integers $\ge 2$. Then the (indecomposable)
rack $\Oc^n_x$ is a disjoint union of $e$ racks isomorphic to
$\Oc^d_x$; in other words, $\Oc^n_x$ is an extension of $\Oc^e_x$
by $\Oc^d_x$ (and vice versa), see \cite[Section 2]{AG1}. Thus,
there is an epimorphism of braided vector spaces $M(\Oc^n_{x},
\sgn) \to M(\Oc^e_{x}, \sgn)$, as well as an  inclusion
$M(\Oc^d_{x}, \sgn) \to M(\Oc^n_{x}, \sgn)$. The techniques
available today do not allow to compute the Nichols algebra
$\toba(\Oc^n_{x}, \sgn)$ from the knowledge of the Nichols algebra
$\toba(\Oc^e_{x}, \sgn)$.
\end{obs}

\begin{obs} In Theorem \ref{dn} we do \emph{not}
claim that the conditions are sufficient. See Tables
\ref{tabladnimpar}, \ref{tabladnpar}. For instance, it is known
that $\dim \toba(\Oc^n_{x}, \sgn) < \infty$ when $n = 3$ -- see
\cite{ms}; for other odd $n$, this is open.
\end{obs}

\begin{table}[t]
\begin{center}
\begin{tabular}{|p{5cm}|p{2cm}|p{1,5cm}|p{1,7cm}|}\cline{1-2}
\hline {\bf Orbit} &    {\bf Isotropy \newline group} & {\bf Rep.}
& $\dim \toba(V)$

\\ \hline  $e$   & $\dn$ & any & $\infty$
\\ \hline  $\Oc_{y^{h}}= \{y^{\pm h}\}$, $h\neq 0$,
\newline $ \mid \Oc_{y^{h}} \mid =2$   & $\Z_n \simeq \langle
y \rangle$ & any &$\infty$

\\  \hline  $\Oc_{x}= \{xy^{h}:0\leq h\leq n-1\}$,
\newline $ \mid \Oc_{x} \mid =n$   & $\Z_2 \simeq \langle
x \rangle$ & $\varepsilon$ &$\infty$

\\ \cline{3-4}   &  & $\sgn$ &  negative\newline braiding
\\ \hline
\end{tabular}
\end{center}

\caption{Nichols algebras of irreducible Yetter-Drinfeld modules
over $\dn$, $n$ odd.}\label{tabladnimpar}

\end{table}

Let us now proceed with the proof of Theorem \ref{dn}.

\pf If $s=\id$, then $q_{ss}=1$ and $\dim \toba(\Oc,\rho)=\infty$,
from Lemma \ref{trivialbraiding}.

We consider now two cases.

\emph{CASE 1:} $n$ odd.

(I) If $s=y^h$, with $1\leq h \leq n$, it is easy to see that
$\Oc_{y^h}=\{ y^h,y^{-h} \}$ and $\dn^{y^h}=\langle y \rangle
\simeq \Z_n$. Then $\widehat{\Z_n}=\{ \chi_l \}_{l=1}^n$, where
$\chi_l(y)=\omega^l$, with $\omega=\exp(\frac{i 2\pi}{n}) \in
\G_{n}$ a primitive $n$-th root of $1$. Let us consider
$M(\Oc_{y^h},\chi_l)$; it is a braided vector space of diagonal
type. If $q_{ss}\neq -1$, then $\dim
\toba(\Oc_{y^h},\chi_l)=\infty$, from Lemma \ref{odd}. Assume
$q_{ss}=-1$; so we have $ -1=\chi_l(s)=\chi_l(y^h)=\omega^{lh}. $
This is a contradiction because $n$ is odd.

(II) If $s=x$, then $\Oc_{x}=\{ x,xy,\dots,xy^{n-1}\}$ and
$\dn^{x}=\langle x \rangle \simeq \Z_2$. Clearly, $\dim
\toba(\Oc_{x},\varepsilon)=\infty$. On the other hand,
$M(\Oc_{x},\sgn)$ is a negative braided vector space, since every
abelian subrack of $\Oc_{x}$ has one element; indeed
$xy^{j}xy^{k}=xy^{k}xy^{j}$, $0\leq j,k \leq n-1$, if and only if
$j=k$.

Therefore, the part (a) of the Theorem is proved.

\begin{table}[t]
\begin{center}
\begin{tabular}{|p{5,4cm}|p{2,1cm}|p{2cm}|p{1,7cm}|}
\hline {\bf Orbit} &    {\bf Isotropy \newline group} & {\bf Rep.}
& $\dim \toba(V)$

\\ \hline  $e$   & $\dn$ & any & $\infty$
\\ \hline  $\Oc_{y^{m}}= \{y^{m}\}$,
 $ \mid \Oc_{y^{m}} \mid = 1$   & $\dn$ & \vspace{1pt}$(V,\rho)\in \widehat{\dn}$,
$\rho(y^{m}) = 1$ & $\infty$

\\ \cline{3-4}     &  & \vspace{0pt} $(V,\rho)\in \widehat{\dn}$,
$\rho(y^{m}) = -1$ & \vspace{0pt} $2^{\dim V}$

\\ \hline  $\Oc_{y^{h}}= \{y^{\pm h}\}$, $h\neq 0, m$,
\newline $ \mid \Oc_{y^{h}} \mid =2$   & $\Z_n \simeq \langle
y \rangle$ & $\chi^j$, \newline $\omega^{hj} = -1$ & 4

\\ \cline{3-4}     &  & $\chi^j$, \newline$\omega^{hj} \neq -1$ & $\infty$

\\ \hline  $\Oc_{x}= \{xy^{2h}:0\leq h\leq m-1\}$
\newline $ \mid \Oc_{x} \mid = m$   & $\Z_2 \times \Z_2 \simeq \newline \langle
x\rangle \oplus \langle y^m \rangle$ & $\varepsilon \otimes
\varepsilon$, \newline $\varepsilon \otimes \sgn$ &$\infty$

\\ \cline{3-4}     &  & $\sgn \otimes \sgn$, \newline $\sgn \otimes \varepsilon$ & negative\newline braiding
\\ \hline  $\Oc_{xy}= \{xy^{2h +1}:0\leq h\leq m-1\}$
\newline $ \mid \Oc_{x} \mid = m$   & $\Z_2 \times \Z_2 \simeq\newline \langle
xy \rangle \oplus \langle y^m \rangle$ & $\varepsilon \otimes
\varepsilon$,
\newline $\varepsilon \otimes \sgn$ &$\infty$

\\ \cline{3-4}     &  & $\sgn \otimes \sgn$, \newline $\sgn \otimes \varepsilon$ & negative\newline braiding
\\ \hline
\end{tabular}
\end{center}

\caption{Nichols algebras of irreducible Yetter-Drinfeld modules
over $\dn$, $n=2m$ even.}\label{tabladnpar}

\end{table}

\medbreak \emph{CASE 2:} $n$ even. Let us say $n=2m$.

(I) If $s=y^m$, then $\Oc_{y^m}=\{y^m\}$ and $\dn^{y^m}=\dn$.
Clearly, $\dim \toba(\Oc_{y^m},\rho)=\infty$, for every $\rho \in
\widehat{\dn}$ with $\rho(s)=\Id$. On the other hand, if $(\rho,V)
\in \widehat{\dn}$ is such that $\rho(s)=-\Id$, then it is
straightforward to prove that $\toba(\Oc_{y^m},
\rho)=\bigwedge(V)$, the exterior algebra of $V$; hence $\dim
\toba(\Oc_{y^m},\rho)=2^{\dim V}$.

(II) If $s=y^h$, $h \neq 0,m$; then $\Oc_{y^h}=\{ y^h ,y^{-h}\}$
and $\dn^{y^h}=\langle y  \rangle \simeq \Z_n$. From Lemma
\ref{odd}, it is clear that $\dim \toba(\Oc_{y^h},\chi_l)=\infty$,
for every $l$ such that $\chi_l(y^h)\neq -1$, i.e.
$\omega^{hl}\neq -1$. On the other hand, it is easy to see that
$\toba(\Oc_{y^h},\chi_l)=\bigwedge(M(\Oc_{y^h},\chi_l))$, hence
$\dim \toba(\Oc_{y^h},\chi_l)=4$, for every $\chi_l$ with
$\chi_l(y^h)=-1$.

(III) If $s=x$, then $\Oc_{x}=\{ x y^{2h} \, : \, 0\leq h \leq m-1
\}$ and $\dn^{x}=\langle x  \rangle \oplus \langle y^m \rangle
\simeq \Z_2 \times \Z_2$. From Lemma \ref{trivialbraiding}, $\dim
\toba(\Oc_{x},\varepsilon \otimes \varepsilon )=\dim
\toba(\Oc_{x},\varepsilon \otimes \sgn)=\infty$.

For the cases $\rho= \sgn \otimes \varepsilon$ or $\sgn \otimes
\sgn$, we note the following fact.
\begin{itemize}
\item[(i)] If $m$ is odd and $0 \leq j,k \leq m-1$, we have that
$$
x y^{2j} x y^{2k} = xy^{2k} xy^{2j} \quad \text{ if and only if }
\quad j=k.$$
\item[(ii)] If $m$ is even and $0 \leq j \leq k \leq m-1$, we have
that
$$
xy^{2j} xy^{2k} = xy^{2k} xy^{2j} \quad \text{ if and only if }
\quad \text{$k=j$ or $j+\frac{m}{2}$}.$$
\end{itemize}
The cases (i) and (ii) say that every maximal abelian subrack of
$\Oc_{x}$ has one and two elements, respectively. Hence, in both
cases the braiding is negative. Indeed, the result is obvious for
the case (i), while in the case (ii) we have that if $t_j:=
xy^{2j}= xy^j \, x \,  (x y^j)^{-1}$ and $t_k:=xy^{2k}=xy^k \, x
\, (xy^k)^{-1}$ commute in $\Oc_{x}$, then $q_{jj}=-1=q_{kk}$ and
$q_{jk}q_{kj}=1$; thus the braiding is negative.

(IV) If $s=xy$, then $\Oc_{xy}=\{ xy^{2h+1} \, : \, 0\leq h \leq
m-1 \}$ and $\dn^{xy}=\langle xy \rangle \oplus \langle y^m
\rangle \simeq \Z_2 \times \Z_2$. The result follows as in (III)
using the isomorphism $\dn \to \dn$, $x\mapsto xy$, $y\mapsto y$.
\epf

\section{On Nichols algebras over semisimple Hopf algebras}\label{nichols-other}

Let $A$ be a Hopf algebra. Let $J\in A\otimes A$ be a twist and
let $A^J$ be the corresponding twisted Hopf algebra. If $A$ is a
Hopf subalgebra of a Hopf algebra $H$, then $J$ is a twist for $H$
and $A^J$ is a Hopf subalgebra of $H^J$. Now, if $A$ is
semisimple, then this induces a bijection
\begin{multline}\label{isom-twist}
\{\text{isoclasses of Hopf algebras with coradical $\simeq A$}\}
\\\overset{\thicksim}{\longrightarrow} \{\text{isoclasses of Hopf
algebras with coradical $\simeq A^J$}\},
\end{multline}
that preserves standard invariants like dimension,
Gelfand-Kirillov dimension, etc. Let now $A = \ku \mathbb A_5$ and
let $J\in A\otimes A$ be the non-trivial twist defined in
\cite{Ni}. By \eqref{isom-twist}, we conclude immediately from
Theorem \ref{mainteor}.

\begin{theorem}\label{twist-a-cinco} Let $H$ be a finite-dimensional Hopf algebra
with coradical isomorphic to $(\ku \mathbb A_5)^J$. Then $H\simeq
(\ku \mathbb A_5)^J$. \qed
\end{theorem}

Again, this is the first classification result we are aware of,
for finite-dimensional Hopf algebras with coradical isomorphic to
a fixed non-trivial semisimple Hopf algebra. Recently, two
semisimple Hopf algebras $B\simeq (\ku\,\mathbb D_3 \times \mathbb
D_3)^{J'}$ and $C\simeq (\ku\,\mathbb D_3 \times \mathbb
D_5)^{J''}$ were discovered in \cite{GN}. Both $B$ and $C$ are
simple, that is they have no non-trivial normal Hopf subalgebra.
Since there are finite-dimensional non-semisimple pointed Hopf
algebras with group either $\mathbb D_3$ or $\mathbb D_5$, there
are finite-dimensional non-semisimple Hopf algebras with coradical
isomorphic either to $B$ or to $C$.

\subsection*{Acknowledgement} We are grateful to Professor John Stembridge
for information on Coxeter groups, in particular reference
\cite{BG}. We thank Mat\'\i as Gra\~na, Sebasti\'an Freyre and
Leandro Vendram\'\i n for interesting discussions.

\end{document}